\documentclass{article}

\usepackage{amssymb}
\newtheorem{theorem}{Theorem}[section]
\newtheorem{proposition}[theorem]{Proposition}
\newtheorem{lemma}[theorem]{Lemma}
\newtheorem{definition}[theorem]{Definition}

\newcommand{\qed}{{\hspace*{\fill}\rule{2.2mm}{2.2mm}}}
\newenvironment{proof}
{\begin{trivlist} \item[] {\bf Proof.}}{\qed \end{trivlist}}
\newcommand{\NN}{{\mathbb N}}
\newcommand{\var}{{\rm var}}
\newcommand{\RR}{{\mathbb R}}

\newcommand{\QQ}{{\mathbb Q}}
\newcommand{\x}{{\bf x}}

\begin{document}

\title{Rational {\L}ukasiewicz logic and DMV-algebras}
\author{\textsc{Brunella Gerla}\cr
{\footnotesize Department of Theorical and Applied Sciences, University of Insubria}\\
{\footnotesize Via Mazzini 5, 21100, Varese, Italy}\\
{\footnotesize {\tt brunella.gerla@uninsubria.it}}\\ \\
}
\date{}
\maketitle
\begin{abstract}
In this paper we present some results on
the variety of divisible MV-algebras.
Any free divisible
MV-algebra is an algebra of continuous piecewise linear functions
with rational coefficients.
Correspondingly, Rational
{\L}ukasiewicz logic is defined and its tautology problem
is shown to be co-NP-complete.
\end{abstract}

\section{Introduction}
Recently,
many-valued logics have been proposed as the mathematical background
for fuzzy logic, and among them a particular importance
has been given to logics in which the conjunction is
modeled by a continuous t-norm.
The best known cases are the minimum, the product and the
 {\L}ukasiewicz t-norm. In particular,
{\L}ukasiewicz logic and its algebraic counterpart,
MV-algebras, have deep relations
with a wide spectrum of mathematical fields, like error-correcting
codes, lattice-ordered group theory and
$C^*$-algebras (see \cite{cigdotmun}).


By McNaughton theorem \cite{McNaughton},
the functions
associated with formulas of {\L}u\-ka\-sie\-wi\-cz logic
are the totality of continuous, piecewise linear functions in which
every piece has integer coefficients.
This fact can be used
to give a formal description of fuzzy rules for fuzzy
control \cite{amapor}. In order to have a more flexible tool,
it seems necessary
to weaken the restriction of integer coefficients,
and consider instead
rational coefficients.

To this purpose
different approaches have been proposed in the literature.
The authors of \cite{agumun01}
introduced
{\L}ukasiewicz propositional logic with one quantified
propositional variable $\exists\mbox{\rm \L}$.
In \cite{dinlet96} Riesz MV-algebras are defined
as a special class of MV-algebras with a family of unary operators,
and are shown to be the MV-algebraic counterpart
of vector lattices over real numbers.
In \cite{hoehle} root operators (in fact, division) are introduced
and in \cite{baazveith} {\L}ukasiewicz logic plus
root operators is shown to correspond to
continuous piecewise linear functions with rational coefficients
and to have the interpolation
property.

In this paper we collect all these results and we
give an equational definition
of root operators, defining the variety of
{\em DMV-algebras} (divisible MV-algebras).
Such structures maintain some basic properties
of MV-algebras, and are {\em intervals}
 of lattice-ordered vector spaces over the rationals
just as MV-algebras
are intervals of lattice-ordered abelian groups
\cite{cigdotmun}.

We extend to DMV-algebras some results holding for
MV-algebras,
like the representation theorem and the correspondence
with divisible l-groups. We further give a direct proof
that the variety
of DMV-algebras is generated by $[0,1]$.
Rational {\L}ukasiewicz logic is then introduced and
is shown to be an extension of Rational Pavelka logic.
In the  last section the tautology problem for
Rational {\L}ukasiewicz logic is shown to be co-NP-complete.

\section{Basic notions: MV\--al\-ge\-bras}

\begin{definition}[\cite{Chang, cigdotmun}]
An {\it MV-algebra} is a structure
$A=(A,\oplus,\neg,0,1)$ satisfying the
following equations:
$$
\begin{array}{ll}
x\oplus (y\oplus z)=(x\oplus y)\oplus z;&
x\oplus y=y\oplus x; \\
x\oplus 0=x; &
x\oplus 1=1; \\
\neg 0=1; &
\neg 1=0; \\
\neg (\neg x\oplus y)\oplus y =\neg (\neg y\oplus x)\oplus x.&
\end{array}
$$
\end{definition}

Boolean algebras
coincide with MV-algebras sa\-ti\-sfy\-ing
the additional equation $x\oplus x=x$ (idempotency).
Each MV-algebra contains as a subalgebra the two-element boolean algebra
$\{0,1\}$.

In any MV-algebra one defines the operations $\odot$ and $\to$
as follows:
$$
\begin{array}{ll}
x\odot y=\neg (\neg x\oplus \neg y)
& \qquad
x\to y=\neg x\oplus y.
\end{array}
$$
Further any MV-algebra $A$ is equipped with the order relation
$$
x\leq y \quad \mathrm{if \; and \; only\; if\; }\quad \neg x \oplus y=1.
$$
Then $A$ becomes a distributive lattice, and
\begin{eqnarray*}
x \wedge y =\inf \{x,y\}= \neg (\neg x\odot y) \odot y;
&\, &
x \vee y =\sup \{x,y\}=\neg (\neg x \wedge \neg y).
\end{eqnarray*}
If its order is total, then the MV-algebra
$A$ is called {\em MV-chain}.
A lattice-ordered group (l-group for short)
$G=(G,0,-,+,\wedge,\vee)$
is an abelian group $(G,0,-,+)$
equipped with a lattice structure $(G,\wedge,\vee)$
such that, for every $a,b,c \in G$,
$c+(a\wedge b)=(c+a)\wedge (c+b)$. An l-group
is said to be totally ordered if the lattice-order
is total.
An element $u\in G$ is a {\em strong unit}
of $G$ if for every
$x\in G$ there exists $n \in \NN$ such that $nu \geq x$.
In \cite{mun_interp} an equivalence functor ${\bf \Gamma}$
from  the category of l-groups with strong unit
to the category of MV-algebras has been constructed.
If $G$ is an l-group and $u$ is a strong unit for $G$,
the MV-algebra
${\bf \Gamma}(G,u)$ has the form
$\{x \in G \mid 0\leq x\leq u\}$ and operations
are defined by $x \oplus y=u \wedge x+y$
and $\neg x=u-x$.
If $A$ is an MV-algebra we shall
denote by $G_A$ the
l-group corresponding to $A$ via ${\bf \Gamma}$.

{\bf Remark.}\,
If $A$ is an MV-algebra,
$x\in A$ and $n\in \NN$ we denote by $n.x$
the element of $A$ inductively defined by $0.x=0$, $(n-1).x=n.x\oplus x$.
Further, we denote by $nx$ the element of $G_A$
defined by $0x=0$ and $(n-1)x=nx+x$.
If $u$ is a strong unit of $G_A$
such that ${\bf \Gamma}(G,u)=A$, it follows that
$n.x=nx\wedge u$.

\medskip
Chang's algebraic completeness theorem \cite{Chang2}
states that
every MV-algebra is a subdirect product of linear MV-chains
and that
an equation in the language of MV-algebras holds
in every MV-algebras if and only if holds in $\QQ\cap [0,1]$.

In \cite{dinola} the author shows that
every MV-algebra is an algebra of functions
taking values in an ultrapower of the interval $[0,1]$.

A standard reference for MV-algebras is \cite{cigdotmun}.

\section{DMV-algebras}
Both the completeness
theorem in \cite{Chang2}
and the representation theorem
in \cite{dinola} are based on results for
the theory of l-groups. In particular the authors
use the result that
every totally ordered group can be embedded in a
totally ordered divisible group, and that quantifier elimination
holds for totally ordered divisible groups.

In this section we shall study algebraic structures
that are more directly connected with divisible groups.
\begin{definition}[Mundici]
A {\em DMV-algebra} $A=(A,\oplus, \neg ,
\{\delta_n\}_{n\in N},0,1)$
is an algebraic structure such that $A^*=(A, \oplus,\neg,0,1)$ is
an MV-algebra and the following hold for every $x\in A$ and $n \in \NN$:
\begin{itemize}
\item[(D1n)]
$n.\delta_n x =x$
\item[(D2n)]
$\delta_n x \odot (n-1).\delta_n x=0$
\end{itemize}
\end{definition}

If $A$ is a DMV-algebra, then the
MV-algebra $A^*$ satisfies the condition of divisibility, i.e.,
for every $n\in N$ and for every $x \in A$ there exists
$y \in A$ such that $n.y=x$.
The MV-algebra $A^*$
is the {\em MV-reduct} of the DMV-algebra $A$.
On the other hand, if $B$ is a divisible MV-algebra
then
by $(B,\delta_n)$ we shall denote
the DMV-algebra obtained from $B$ by the introduction
of the new connective $\delta_n$ for every $n>0$.

{\bf Example}\,
For each $k=1,2,\ldots,$ the set
$$
{\mbox {\L}}_{k+1}=\{0,\frac 1{k},\ldots,\frac {k-1}{k},1\},
$$
equipped with the operations
$$
x\oplus y =\min \{1,x+y\},\;\;\; x\odot y=\max
\{0,x+y-1\},\;\; \neg x=1-x
$$
is a linearly ordered MV-algebra (also called \textit{MV-chain}),
but cannot be enriched to a DMV-algebra.
The set of all rationals
between $0$ and $1$
where each $\delta_n$ is interpreted as division
by $n$,  is a DMV-algebra that we shall denote by
$({\bf \Gamma}(\QQ,1),\delta_n)$.
In this case, Axioms $(D1n)$ and $(D2n)$
state that the sum of $n$ copies
of $x/n$ coincides with $x$.

\begin{proposition}
\label{unique}
Let $A$ be a DMV-algebra,
let $A^*$
be its MV-reduct and let
$(G,u)$ be
the unique l-group with strong unit $u$
such that ${\bf \Gamma}(G,u)=A^*$.
Then, for every $x\in A$,
\begin{itemize}
\item[(i)]
$m\delta_m x=x$
\item[(ii)]
$\delta_m x$ is the unique element of $A$
satisfying axioms (D1) and (D2).
\end{itemize}
\end{proposition}
\begin{proof}
\begin{itemize}
\item[(i)]
If $x\in A$, equations (D1) and (D2) become
\begin{eqnarray}
\label{def1}&&m\delta_m x\wedge u=x\\
\label{def2}&&(\delta_m x + ((m-1)\delta_m x\wedge u)-u) \vee 0 =0
\end{eqnarray}
whence $\delta_mx +((m-1)\delta_m x\wedge u)\leq u$
and, by definition of l-group,
$m\delta_m x \wedge (u+\delta_m x)\leq u$.
Since $\delta_n x \geq 0$, then
$u+\delta_n x\geq u$ and
$$
m\delta_m x \wedge u \leq m\delta_m x \wedge( u+\delta_m x)
\leq u
$$
and hence, from (\ref{def1}),  $m\delta_m x =x$.
\item[(ii)]
For every $y\in A\subseteq G$ satisfying (D1) and (D2),
$my\wedge u=x$ and
$(y+((m-1)y\wedge u)-u)\vee 0=0$.
Repeating the same argument as above,
$my=x=m\delta_m x$ and then $y=\delta_m x$.
\end{itemize}
\end{proof}

Chang's distance function $d: A\times A \to A$ is defined by
$$
d(x,y)= (x \odot \neg y) \oplus (y \odot \neg x).
$$

\begin{proposition}
\label{delta-dist}
Let $A$ be a DMV-algebra and let $x,y \in A$.
\begin{itemize}
\item[(i)]
If $x \odot y=0$ then $\delta_n(x\oplus y)=\delta_n x \oplus \delta_n y$

\item[(ii)]
$\delta_n d(x,y)=d(\delta_n x, \delta_n y)$

\end{itemize}

\end{proposition}
\begin{proof}
We shall give the proof for
the case $n=2$.
This can be generalized to
every $n>0$.
\begin{itemize}
\item[(i)]
If $x,y\in A$,
$$
(\delta_2 x\oplus \delta_2 y)\oplus  (\delta_2 x\oplus \delta_2 y)=
(\delta_2 x\oplus \delta_2 x) \oplus   (\delta_2 y\oplus \delta_2 y)=
x \oplus y.
$$
Further note that in every MV-algebra, if $a\odot b=0, a\odot a=0$
and $b\odot b=0$ then
$(a\oplus b)\odot (a\oplus b)=(a \oplus a) \odot (b\oplus b)$. Thus,
since $\delta_2 x \odot \delta_2 y\leq x \odot y=0$,
$$
(\delta_2 x\oplus \delta_2 y)\odot (\delta_2 x\oplus \delta_2 y)=
x \odot y=0.
$$
Therefore, $\delta_2(x\oplus y)=\delta_2 x \oplus \delta_2 y$
\item[(ii)]
Let $a,b$ elements of $[0,1]$ such that
$a \odot b=0$. In $[0,1]$ we have $d(a,b)=|a-b|$ and hence
$$
d(a,b)\odot d(a,b)=|a-b|\odot |a-b|= 2|a-b|-1 \vee 0.
$$
If $a\leq b$ then $2|a-b|-1 \vee 0=2(a-b)-1 \vee 0 = 2a-2b-1 \vee 0$ and
since $2a-1=0$ then $d(a,b)\odot d(a,b)=0$. Analogously
the same conclusion can be drawn in case $a\geq b$.

By the Chang representation theorem it follows that
if $\delta_2 x \odot \delta_2 x =0$,
$$
d(\delta_2 x, \delta_2 y)\odot d(\delta_2 x, \delta_2 y)=0.
$$
Further,
$$
d(\delta_2 x, \delta_2 y)\oplus  d(\delta_2 x, \delta_2 y)=
d(2 \delta_2 x, 2\delta_2 y)=d(x,y)
$$
whence the claim follows.
\end{itemize}

\end{proof}

\begin{definition}
If $A$ and $B$ are DMV-algebras, a function
$f:A\to B$ is a {\em homomorphism} of DMV-algebras if
$f$ is a MV-homomorphism from $A^*$
to $B^*$ and
for every $x \in A$,
$$
f(\delta_n x)=\delta_n f(x)
$$
\end{definition}

\begin{definition}
A subset $J$ of a DMV-algebra $A$
is said to be a {\em DMV-ideal} of $A$
if it is an ideal of the MV-reduct $A^*$, that is:
\begin{itemize}
\item
$0 \in J$
\item
For every $x \in J$ and $y\leq x$ then $y \in J$
\item
If $x,y \in J$ then $x\oplus y \in J$
\end{itemize}
\end{definition}

Note that if $J$ is an ideal and $x\in J$ then also $\delta_n x\in J$
for every $n$.
A DMV-ideal $J$ is a {\em prime ideal} iff it is not trivial and
for every $x,y \in A$, either $x \odot \neg y\in J$
or $y \odot \neg x \in J$.

\begin{proposition}
\label{congr}
Let $I$ be an ideal of $A$.
The binary relation $\equiv_I$ on $A$ defined by
$x \equiv_I y$ if and only if $d(x,y)\in I$ is a congruence
relation.
\end{proposition}

\begin{proof}
Indeed $\equiv_I$ is a congruence on the MV-reduct
$A^*$.
Further, if $x,y \in A$ and $x\equiv_I y$ then,
by Proposition \ref{delta-dist},
$d(\delta_n x, \delta_n y)=\delta_n d(x,y) \leq d(x,y)$,
hence $d(\delta_n x, \delta_n y)\in I$ and $\delta_n x=\delta_n y$.
\end{proof}

Let $I$ be a DMV-ideal of $A$ and
let $\pi:x\in  A \mapsto [x]_I \in A/I$.
Then $Ker(\pi)=\{x\in A \mid [x]_I=[0]_I\}=
\{x\in A\mid d(x,0)\in I\}=I$.
Vice-versa, if $f$ is a DMV-homomorphism,
then $Ker(f)=\{x\in A\mid f(x)=0\}$ is
a DMV-ideal. We get
\begin{proposition}
$I$ is a DMV-ideal of $A$ if and only if
there exists a DMV-homomorphism $f$
such that $I=Ker(f)$.
\end{proposition}

By Proposition \ref{congr}, if $I$ is an DMV-ideal of
$A$, then
setting
\begin{itemize}
\item
$[x]_I\oplus [y]_I=[x\oplus y]_I$
\item
$\neg[x]_I=[\neg x]_I$
\item
$\delta_n[x]_I=[\delta_n x]_I,$
\end{itemize}
the structure $(A/I, \oplus,\neg, \{\delta_n\}, [0]_I)$
is a DMV-algebra.
Further, the quotient $A/I$ is totally ordered
iff $I$ is a prime ideal.
The proof of the following Proposition is the same as
for MV-algebras:
\begin{proposition}
\label{primeideal}
Let $A$ be a DMV-algebra and $I$ an ideal of $A$.
If $z\not\in I$ then
there exists a prime ideal $P$ of $A$ such that
$I\subseteq P$ and $z\not\in P$.
\end{proposition}

Then we can extend to DMV-algebras
the Chang Representation theorem.
\begin{theorem}
\label{linear}
Every DMV-algebra is the subdirect product of linear DMV-algebras.
\end{theorem}

The functor ${\bf \Gamma}$
induces a correspondence
between DMV-algebras and divisible l-groups:

\begin{definition}[\cite{mun_interp}]
A {\em good sequence} of a DMV-algebra $A$ is a sequence
$(a_1,$ $\ldots,$ $a_n)$ of elements of $A$ such that
for every $i=1,\ldots,n-1$, $a_i\oplus a_{i+1}=a_i$.
\end{definition}
If $A$ is linear then every good sequence
has the form $(1^p,a)=
(\underbrace{1,\ldots,1}_{p\hbox{ times}},a)$ with $a\in A$.
Further, if $(a_1,\ldots,a_n)$ is a good sequence of $A$
then also $(a_1,\ldots,a_n,0)$ is a good sequence.

\begin{proposition}
\label{lindiv}
Let $A$ be a totally ordered DMV-algebra. Then there exists a
totally ordered divisible group $G$ together with a strong unit $u$
such that $A=\{x\in G \mid 0\leq x\leq u\}$.
\end{proposition}
\begin{proof}
If $A$ is a totally ordered  DMV-algebra,
then thge MV-reduct $A^*$ is a MV-chain
and $A^*$ is isomorphic with ${\bf \Gamma}(G_{A^*},(1,0))$
(see \cite{mun_interp}).
Let $u=(1,0)$.
For every $n \in N$ and for every $x \in [0,u]$ there exists
$y\in [0,u]$ such that $ny=x$.
Since $u$ is a strong unit and $G$ is linear,
for every $x \in G$ there
exists an integer number $n_x$ such that
$n_xu \leq x< (n_x +1)u$. Let $x'=x-n_xu \in [0,u]$. Then
let $y'$ such that $ny'=x'$ and let $u'$ such that
$nu'=u$. Then the element $n_xu'+y'$ is such that
$n(n_xu'+y')=x$, hence the totally ordered
group $G_{A^*}$ is divisible.
\end{proof}

\begin{theorem}
Let $A$ be a DMV-algebra. Then there exists a
unique divisible l-group $G$
together with a strong unit $u$ for $G$ such that
$A=\{x\in G \mid 0\leq x\leq u\}$.
\end{theorem}
\begin{proof}
From Theorem \ref{linear}, $A$ is a subdirect product
of totally ordered DMV-algebras $(A_i)_{i\in I}$
and every $A_i$ is equal to ${\bf \Gamma}(G_i,u_i)$
with $G_i$ totally ordered divisible group
(Proposition \ref{lindiv}), hence
\begin{equation}
\label{subdirect}
A\subseteq \prod_{j\in J}A_j \subseteq \prod_{j\in J}G_j.
\end{equation}
Let $u=(u_j)_{j\in J}$.
By \cite{mun_interp}, if $G$ is the group
generated by $A$ in $\prod_{j\in J}G_j$ and if
$$
G^+=\left\{ a_1+\ldots + a_n \mid
(a_1,\ldots,a_n) \hbox{ good sequence of } A \right\}
$$
then $G=G^+-G^+$, $u$ is a strong unit of $G$ and
${\bf \Gamma}(G,u)=A^*$.
There remains to show that $G$ is divisible, i.e.
if $m\in \NN$, for every $x\in G$
there exists $y\in G$ such that $my=x$.
It is enough to restrict to
$x\in G^+$.

Let  $a_1,\ldots,a_r$ be a
good sequence of $A$ ($a_r\neq 0$) and
$x=a_1+\ldots +a_r \in G^+$.
Then $\delta_n a_i\in A\subseteq G$ for every $i=1,\ldots,r$.
Let
$$
y=\delta_n a_1 +\ldots +\delta_n a_r \in G.
$$
By Proposition \ref{unique}(i) we have $ny=x$.
\end{proof}

Since for divisible, totally ordered abelian groups
the quantifier elimination theorem holds, then
a universal sentence $\chi$ is satisfied by every
divisible totally ordered abelian group if and only if it is
satisfied by $\QQ$.


\begin{definition}
\label{DMV-equation}
A DMV-equation in the variables $X_1,\ldots,X_n$
is an expression $\tau=\sigma$, where
$\tau$ and $\sigma$ are terms
over the alphabet $\{\oplus, \neg, 0,1\}\cup \{\delta_n\}_{n\in \NN}$
with variables among $X_1,\ldots,X_n$ ({\em DMV-terms}).
A DMV-equation $\tau=\sigma$ is {\em satisfied} by a DMV-algebra $A$
if, for every $n$-tuple $(a_1,\ldots,a_n)\in A^n$,
$\tau(a_1,\ldots,a_n)=\sigma(a_1,\ldots,a_n)$, where
$\tau(a_1,\ldots,a_n)$ and $\sigma(a_1,\ldots,a_n)$
are elements
of $A$ obtained by substituting
$X_1,\ldots,X_n$ by $a_1,\ldots, a_n$ in $\tau$ and $\sigma$.
\end{definition}
Repeating the same argument as for MV-algebras
(\cite{Chang2}), we have

\begin{theorem}
\label{represent}
A DMV-equation is satisfied by every MV-algebra if and only if
it is satisfied by the DMV-algebra $({\bf \Gamma}(\QQ,1),\delta_n)$.
\end{theorem}

\subsection{Varieties and quasi-varieties of DMV-algebras}
\label{variety}
Since DMV-algebras have an equational definition,
the class of all DMV-algebras is a variety.
By Theorem \ref{represent} we have
\begin{theorem}
The variety of DMV-algebras is generated by $[0,1]\cap \QQ$.
\end{theorem}

It is possible to give an alternative proof of this theorem,
by translating
any equation of DMV-algebras in a quasi-equation of MV-algebras.

Indeed, suppose that $\tau=1$ is a DMV-equation and let
${\mathcal T}$ be the parsing tree of $\tau$,
that is, ${\mathcal T}$ is a tree
which nodes are subformulas of $\tau$ and such that
each node has as children its direct subformulas.
Leaves of ${\mathcal T}$
are all the occurrences of variables occurring in $\tau$.

Let us display the occurrences
of variables $x_1,\ldots, x_n$ in any formula $\varphi$ by
writing $\varphi(x_1,\ldots,x_n)$.
Suppose that
$\delta_{i_1}\tau_{i_1}, \ldots, \delta_{i_m}\tau_{i_m}$
is an enumeration of all nodes of $T$ that begin with a symbol
$\delta$. Each of these nodes $\delta_{i_j}\tau_{i_j}$
has a unique child $\tau_{i_j}$.
Let us introduce $m$ new variables in order to eliminate
the occurrences of $\delta$:
if $\tau=\tau(x_1,\ldots,x_n)$, let
$\tau^*(x_1,\ldots, x_n,z_1,\ldots,z_m)$ be the formula
obtained by substituting every subformula $\delta_{i_j}\tau_{i_j}$
by $z_j$.

If $\tau_{i_j}$ is disjoint from any other subformula
in $\{\tau_{i_1},\ldots, \tau_{i_m}\}\setminus \{\tau_{i_j}\}$
then we denote by $\sigma_1(z_j)$
the MV-equation $\{i_j.z_j=\tau_{i_j} \}$
and by $\sigma_2(z_j)$ the MV-equation $\{ z_j \odot (i_j-1)z_j=0 \}$.
Otherwise, suppose that there exists $h_1,\ldots, h_l
\in \{i_1,\ldots,i_m\}$ such that $\tau_{h_1},\ldots, \tau_{h_n}$
are subformulas of $\tau_{i_j}$. By induction,
let $\tau^*_{i_j}$ be obtained by substituting
each $\tau_{h_k}$ by $z_{h_k}$ and let
$\sigma_1(\tau_{i_j})$ be the MV-equation
$\{ i_j.z_{i_j}=\tau^*_{i_j} \}$ and
$\sigma_2(\tau_{i_j})$ be the MV-equation
$\{ z_{i_j} \odot (i_j-1).z_{i_j}=0 \}$.

The equation $\tau=1$ holds in a DMV-algebra $A$ if and only if
the quasi-equation
\begin{equation}
\label{quasi-eq}
\begin{array}{ll}
\hbox{ IF }
\left.
\begin{array}{c}
\sigma_1{\tau_{i_1}}\hspace{-.2cm}=1
\hbox{ AND } \sigma_2{\tau_{i_1}}\hspace{-.2cm}=1\\
\vdots \\
\hbox{ AND } \sigma_1{\tau_{i_m}}\hspace{-.2cm}=1\,
\hbox{ AND } \sigma_2{\tau_{i_m}}\hspace{-.2cm}=1\,
\end{array}
\right\}
&
\hbox{ THEN }
\, \tau^*=1
\end{array}
\end{equation}
holds in the MV-reduct $A^*$. Since the quasi-varieties of MV-algebras
is generated by $\QQ\cap [0,1]$, quasi-equation (\ref{quasi-eq})
fails in an MV-algebra if and only if it fails in $\QQ\cap [0,1]$.

\bigskip
In \cite{dinola} it is shown that every MV-algebra is
an algebra of functions over an ultrapower of $[0,1]$. This
is equivalent to saying that the quasi-variety generated by $[0,1]$
is the whole variety of MV-algebras. The proof of this theorem
can be adapted to DMV-algebras in
the following way:

Let $A$ be a DMV-algebra. Then $A$ is a subdirect product of
totally ordered DMV-algebras
$A_i={\bf \Gamma}(G_i,u_i)$.
Since each $G_i$ is
a totally ordered divisible group, then
it is elementarily equivalent to the additive group $\RR$ of real numbers
with natural order. Then ${\bf \Gamma}(G_i,u_i)$ is elementarily equivalent
to the MV-algebra $[0,1]$ and hence, by Frayne's theorem
(see for example \cite{chakie}),
it is elementarily embeddable in a suitable ultrapower $[0,1]^{*_i}$
of $[0,1]$. Therefore, since
$$
A\subseteq \prod_{j\in J}A_j \subseteq \prod_{j\in J}{\bf \Gamma}(G_j,u_j)
\subseteq \prod_{j\in J} [0,1]^{*_j}
$$
and applying the joint embedding property of first-order logic,
there exists an ultrapower of $[0,1]^*$ of $[0,1]$
only depending on $A$ such that $A\subseteq [0,1]^*$.

\section{Rational {\L}ukasiewicz logic}
Formulas of Rational {\L}ukasiewicz calculus
are built from
the connectives of negation ($\neg$),
implication ($\to$),
and division ($\delta_n$)
in the usual way.
An axiom is a formula that can be written in any one
of the following ways, where $\varphi$, $\psi$ and $\gamma$
denote arbitrary formulas:
\begin{itemize}
\item[A1)]
$\varphi \to(\psi \to \varphi)$
\item[A2)]
$(\varphi \to \psi)\to
((\psi \to \gamma)\to (\varphi \to \gamma))$
\item[A3)]
$((\varphi \to \psi)\to \psi)
\to ((\psi \to \varphi)\to \varphi)$
\item[A4)]
$(\neg\varphi \to \neg\psi)\to (\psi \to \varphi)$
\end{itemize}
plus, writing $\varphi\oplus \psi$ as
an abbreviation of $\neg \varphi \to \psi$,
\begin{itemize}
\item[A5)]
$\underbrace{\delta_n \varphi \oplus \ldots \oplus \delta_n \varphi}_
{n \hbox{ times}}\to \varphi$
\item[A6)]
$\varphi\to
\underbrace{\delta_n \varphi \oplus \ldots \oplus \delta_n \varphi}_
{n \hbox{ times}}$
\item[A7)]
$\neg \delta_n \varphi \oplus \neg
\underbrace{(\delta_n \varphi \oplus \ldots \oplus
\delta_n \varphi)}_{n-1 \hbox{ times}}.$
\end{itemize}
We shall denote by ${\bf 1}$
the formula $X\to (X \to X)$
where the variable $X$ is fixed once and for all.
Proofs and provability are as usual;
if $\Gamma$ is a set of formulas,
$\Gamma \vdash \varphi$ means that $\Gamma$ proves $\varphi$
(or $\varphi$ is provable from $\Gamma$), that is
there exists a sequence of formulas
$\gamma_1,\ldots,\gamma_u$
such that $\gamma_u=\varphi$ and every $\gamma_i$ either is an
axiom of rational {\L}ukasiewicz logic, or belongs to $\Gamma$
or is obtained from $\gamma_{i_1},\gamma_{i_2}$
($i_1,i_2<i$) by modus ponens.
$\varphi$ is {\em provable} ($\vdash \varphi$)
if is provable from
the emptyset.

Let $Form$ be the set of Rational {\L}ukasiewicz formulas
and let $\equiv$ be the binary relation over $Form$
defined by
$\varphi \equiv \psi$ if and only if
$\varphi \to \psi$ and $\psi \to \varphi$ are provable.
Then $\equiv$ is an
equivalence relation and if
$\varphi$ and $\psi$ are provable formulas then
$\varphi\equiv\psi$.
\begin{proposition}
[Lindenbaum algebra]
\label{lind}
The set
${\mathcal L}=Form /\equiv$
equipped with the operations
\begin{eqnarray*}
\neg [\varphi]_\equiv=[\neg \varphi]_\equiv; &\quad
[\varphi]_\equiv\oplus [\psi]_\equiv=
[\varphi \oplus \psi]_\equiv; &\quad
\delta_n [\varphi]_\equiv=[\delta_n \varphi ]_\equiv
\end{eqnarray*}
is a DMV-algebra where
$1=\{[\varphi]_\equiv\mid \varphi \hbox{ is provable}\}
=[{\bf 1}]_\equiv$.
\end{proposition}
\begin{proof}
Since a similar result holds for {\L}ukasiewicz logic,
we have to prove that
${\mathcal L}$ satisfies $D1n$ and $D2n$.
Indeed for every $[\varphi]_\equiv\in {\mathcal L}$,
by Axioms A5 and A6,
$$
n.\delta_n [\varphi]_\equiv= [n.\delta_n \varphi]_\equiv
=[\varphi]_\equiv
$$
and by Axiom A7,
\begin{eqnarray*}
\delta_n [\varphi]_\equiv\odot (n-1). \delta_n [\varphi]_\equiv&=&
\neg(\neg \delta_n [\varphi]_\equiv \oplus \neg
(n-1).\delta_n [\varphi]_\equiv)=\\
&&\neg [\neg \delta_n \varphi \oplus \neg (n-1)\delta_n
\varphi]_\equiv=\neg 1.
\end{eqnarray*}
\end{proof}
Interpretation of connectives of Rational {\L}ukasiewicz logic
is given by
\begin{definition}
\label{definition:assignment}
An {\em assignment} is a function
$v: Form \to [0,1]$ such that
\begin{itemize}
\item
$v(\neg \varphi)=1- v(\varphi)$
\item
$v(\varphi \to \psi)= \min(1-v(\varphi)+ v(\psi),1)$
\item
$v(\delta_n \varphi)=\frac {v(\varphi)}{n}$.
\end{itemize}
\end{definition}
Every function $\iota$ from the set of variables to $[0,1]$
is uniquely extendible to an assignment $v^\iota$.
For each point $\x=(x_1,\ldots,x_n) \in [0,1]^n$
let $\iota_\x$ be the function mapping each variable $X_j$ into
$x_j$. Fixed $n$ with each formula $\varphi$
with $|\var(\varphi)|\leq n$
it is possible to associate the function
$$
f_\varphi: \x \in [0,1]^n \mapsto v^{\iota_{\x}}(\varphi)\in [0,1]
$$
satisfying the following conditions:
\begin{itemize}
\item $f_{X_i}(x_1,\ldots,x_n) =  x_i$ = the $i$th projection.
\item $f_{\neg\varphi} = 1 - f_\varphi$.
\item $f_{(\varphi \to \psi)} = \min(1,1-f_\varphi + f_\psi)$
\item
$f_{(\delta_n \varphi)}=\frac{f_\varphi}{n}$.
\end{itemize}

A formula $\varphi$ with $|\var(\varphi)|<n$
is {\em satisfiable} iff there exists $\x \in [0,1]^n$ such
that $f_\varphi(\x)=1$. $\varphi$ is a {\em tautology}
iff for every $\x\in [0,1]^n$, $f_\varphi(\x)=1$.
An assignment $v$ is a {\em model} of a set
of formulas $\Gamma$ if for every $\tau \in \Gamma$,
$v(\tau)=1$.

\begin{theorem}[Completeness]
If $\varphi$ is a tautology of Rational {\L}ukasiewicz calculus, then
$\varphi$ is provable.
\end{theorem}
\begin{proof}
Suppose that $\varphi$ is not provable. Then the equation
$\varphi = 1$ is not true in the Lindenbaum DMV-algebra ${\mathcal L}$
of Proposition \ref{lind}, and so
by Theorem \ref{represent}, $\varphi \neq 1$ in
$({\bf \Gamma} (\QQ,1),\delta_n)$. This means that there
exists ${\bf y}\in [0,1]^n$ such that $f_\varphi({\bf y})<1$,
hence $\varphi$ is not a tautology.
\end{proof}

Then, $\varphi \equiv \psi$ if and only if
$\vdash \psi\to \varphi$ and $\vdash \varphi \to \psi$,
if and only if, for every assignment $v$, $v(\varphi)=v(\psi)$,
if and only if $f_\varphi=f_\psi$.
\subsection{Free DMV-algebras}

The Lindenbaum algebra of Proposition \ref{lind}
is the free DMV-algebra $Free_\omega$
over a denumerable set of generators.
In this section we shall describe the free DMV-algebra over a finite
number of generators in terms of continuous piecewise
linear functions.

\medskip

A direct inspection shows that
every function $f_\varphi$ is a continuous
piecewise linear function, where each piece has
rational coefficients.

McNaughton theorem
(\cite{McNaughton} and \cite{mun_constr} for a constructive proof)
states that
a function is associated with a
{\L}ukasiewicz formula
if and only if it is a continuous piecewise linear
function, each piece having integer coefficients.
In \cite{baazveith}, the authors showed that for
every continuous piecewise
linear function $f$ with rational coefficients
there exists a {\L}ukasiewicz formula $\tau$
with division operators  such that $f=f_\tau$.
The proof can be summarized as follows:

Let $f:[0,1]^n \to [0,1]$ be a continuous piecewise linear
function, such that each piece has rational coefficients.
Further, let $s$ be an integer number such that
$s\cdot f:\x \in [0,1]^n \mapsto s\cdot f(\x)\in [0,s]$
is a continuous function with integer coefficients
(for example $s$ is the least common multiple
of the coefficients' denominators of pieces of $f$).

For every $i=0,\ldots,s-1$, let
$$
f_i: \x \in [0,1]^n\mapsto
((s\cdot f(\x)-i)\wedge 1)\vee 0 \in [0,1].
$$
For every $\x \in [0,1]^n$ such that $f(\x) \in [i,i+1]$,
we have $s\cdot f(\x)=i+f_i(x)$.
Since $f_i$ are continuous functions with integer coefficients
there exist MV-terms $\psi_i$ such that
$f_i=f_{\psi_i}$.
If $g:[0,1]^n \to [0,1]$ is any function, let us define
\begin{eqnarray*}
Supp(g)&=&\{\x \in [0,1]^n \mid g(\x)>0\}\\
Supp^{<1}(g)&=&\{\x \in [0,1]^n \mid 0<g(\x)<1\}.
\end{eqnarray*}
We have, for every $i=1,\ldots,s-1$,
$$
Supp^{<1}(f_i)\subseteq Supp(f_i)\subseteq Supp(f_{i-1}).
$$
Indeed
$$
Supp(f_i)=\{\x \in [0,1]^n\mid s\cdot f(\x)>i\}\subseteq
\{\x \in [0,1]^n\mid s\cdot f(\x)>i-1\}.
$$
Further, for any $i\neq j$, $Supp^{<1}(f_i)\cap Supp^{<1}(f_j)=\emptyset$.

\begin{proposition}
\label{eqpiani}
In accordance with the previous notation,
if $f:[0,1]^n \to [0,1]$ is a continuous
piecewise linear function with rational coefficients,
then for every $\x\in [0,1]^n$,
$$
f(\x)=f_\varphi(\x)\qquad \quad \hbox{ where }\quad
\varphi=\bigoplus_{i=1}^{s-1}\delta_s \psi_i.
$$
\end{proposition}
\begin{proof}
Suppose that $\x \in [0,1]^n$ and $f(\x)=0$. Then
for every $i=0,\ldots,s$, $f_i(\x)=0$ whence $f_{\psi_i}=0$
and $f_\varphi=0$.

If $f(\x)=1$ then for every $i=0\ldots,s-1$
$f_{\psi_i}=1$ whence $f_{\delta_s \psi_i}=1/s$ and
$f_\varphi=1$.

Suppose now that there exists $i\in \{0,\ldots,s-1\}$
such that $i<s\cdot f(\x)<i+1$. Then $x\in Supp(f_i)$ and
$f(\x)=f_i(\x)+i/s$.
For every $j\geq i+1$, we have $s\cdot f(\x)-j\leq
s\cdot f(\x)-i-1<0$ whence $f_j(\x)=0$.
Further, for every $j\leq i-1$,
we have $s\cdot f(\x)-j\geq s\cdot f(\x)-i+1>1$
whence $f_j(\x)=1$.

The last case to consider is when $0<s\cdot f(\x)=i<1$. Then
for every $j\leq i$, $f_j(\x)=0$ and for every $k>i$,
$f_k(\x)=1$.

\end{proof}

\begin{theorem}
The free DMV-algebra over $n$ generators is
the algebra of all functions from $[0,1]^n$
to $[0,1]$ that are continuous, piecewise linear and such that
each linear piece has rational coefficients.
\end{theorem}

\begin{proof}
Let $\mathcal{RM}_n$ denote the set of continuous
piecewise linear function with rational coefficients
over $n$ variables and let ${\mathcal X}=\{x_1,\ldots,x_n\}$
the set of variables.
By identifying each variable $x_i$ with the
$i$-th projection, ${\mathcal X}$ is included in $\mathcal{RM}_n$.
If $A$ is any DMV-algebra and $h$ is a map from $\mathcal{X}$ to
$A$, then, for every $f_\varphi\in \mathcal{RM}_n$,
the map
$$
\beta_h(f_\varphi(x_1,\ldots,x_n))=f_\varphi(h(x_1),\ldots,h(x_n))
$$
is a DMV-homomorphism such that
$\beta_h(x_i)=h(x_i)$ for every $x_i \in {\mathcal X}$.
If $\gamma:\mathcal{RM}_n \to A$ is any DMV-homomorphism
such that $\gamma(x_i)=h(x_i)$, then
$$
\begin{array}{ll}
\gamma(f_\varphi(x_1,\ldots,x_n))&=f_\varphi(\gamma(x_1),\ldots,
\gamma(x_n))=\\
&=f_\varphi(h(x_1),\ldots,h(x_n))=\beta_h(f_\varphi(x_1,\ldots,x_n)
\end{array}
$$
hence $\gamma=\beta_h$.
\end{proof}

\subsection{Pavelka-style Completeness}

In \cite{pavelka} the author,
starting from the notion of many valued rules of inference,
defined a class of
complete residuated lattice-valued propositional
calculi and introduced degrees of
provability and degrees of validity.
Then he proved that in {\L}ukasiewicz propositional calculus,
enriched by a denumerable set of rational
constants (what in \cite{Hajek} is called
Rational Pavelka Logic), the degree of provability
of each formula coincides
with the degree of validity
({\em Pavelka-style completeness}).

We shall show that Rational {\L}ukasiewicz logic is a
proper extension of
Rational Pavelka logic. Indeed every formula of
Rational Pavelka logic can be
expressed in Rational {\L}ukasiewicz language,
and, after defining the degree of provability and the degree of
truth,
we shall prove that the completeness with
respect to this degrees
still holds. We shall adapt to our context
the arguments in \cite{cigdotmun, Hajek}.

\begin{definition}
An R{\L}-theory $T$
is a set of Rational {\L}ukasiewicz formulas such that
\begin{itemize}
\item
All axioms belong to $T$;
\item
If $\varphi \to \psi \in T$ and $\varphi \in T$
then $\psi \in T$.
\end{itemize}
\end{definition}

If $T$ is an R{\L}-theory, let us denote by $[T]$ the
set $\{[\varphi]_\equiv \mid \varphi \in T\}$.
Then $T$ is an R{\L}-theory if and only if
$\neg [T]=\{ [\neg \varphi]_\equiv \mid \varphi \in T\}$
is an ideal of the Lindenbaum algebra
in Proposition \ref{lind}.
If $X$ is any set of formulas, then
the {\em R{\L}-theory ${\bf Th}(X)$ generated by $X$}
 is the smallest R{\L}-theory containing $X$.


An R{\L}-theory $T$ is {\em consistent}
if there exists a formula $\varphi$ such that $\varphi \notin T$.
Following \cite{munpan}, an R{\L}-theory is {\em prime}
if it is consistent and for every pair of formulas $\varphi$ and $\psi$,
either $\varphi \to \psi \in T$ or $\psi\to \varphi\in T$.

By Proposition \ref{primeideal},
if $T$ is a consistent R{\L}-theory then there exists
a prime R{\L}-theory $T'$ such that $T'\supseteq T$.

\begin{definition}
Let $\Gamma$ be an R{\L}-theory
and $\varphi$
an Rational {\L}ukasiewicz formula. For every $r/s \in \QQ\cap[0,1]$
we shall denote
by $r/s$ the formula $r. (\delta_s {\bf 1})$. Then,

- the {\em truth degree} of $\varphi$ over $\Gamma$ is
$||\varphi||_\Gamma=\inf \{ v(\varphi)\mid v
\hbox{ is a model of } \Gamma\}$;

- the {\em provability degree} of $\varphi$ over $\Gamma$ is
$|\varphi|_\Gamma=\sup \{ r \mid r\to\varphi \in \Gamma \}$.
\end{definition}
Note that if $\varphi \in \Gamma$ then
by Axiom A1, ${\bf 1} \to \varphi \in \Gamma$.
Hence $|\varphi|_\Gamma=1$.

\medskip
In order to prove the completeness theorem,
we recall the following results
holding for Rational Pavelka logic, that can be
easily generalized for Rational {\L}ukasiewicz logic.
\begin{lemma}
\label{lemmapavelka}
Let $T$ be an R{\L}-theory.
\begin{itemize}
\item[(a)]
If $T$ does not contain $(r \to \varphi)$ then
the R{\L}-theory ${\bf Th}(T\cup \{\varphi \to r\})$
generated by $T \cup \{\varphi \to r\}$ is consistent.
\item[(b)]
If $T$ is prime, for each $\varphi$
$$
|\varphi|_T=\sup \{r \mid r \to \varphi \in T \}=
\inf \{ s \mid \varphi \to s \in T\}.
$$
\end{itemize}
\end{lemma}

\begin{theorem}
\label{evaluation}
If $T$ is a prime R{\L}-theory,
the function $e: \varphi\in Form \to |\varphi|_T\in [0,1]$
is an assignment. That is,
\begin{eqnarray*}
|\neg \varphi|_T=1-|\varphi|_T, &|\varphi\to \psi|_T=
|\varphi|_T\to |\psi|_T,&
|\delta_n \varphi|_T=\frac{|\varphi|_T}{n}
\end{eqnarray*}
hence $e$ is a model of $T$.
\end{theorem}

\begin{proof}
Since the theorem holds for Rational Pavelka
logic, we only have to prove
$|\delta_n \varphi|_T=\frac 1{n}\cdot
|\varphi|_T$.

Since $\vdash (t\to\delta_n \varphi) \to (nt\to\varphi)$,
\begin{eqnarray*}
\frac 1{n}\cdot |\varphi|_T&=&
\frac{\inf \{ s\mid \varphi\to s \in T\}}{n}\\
&=& \inf\{ \frac s{n} \mid  \varphi \to s \in T\}=
\inf\{ t \mid  \varphi \to nt \in T\} \leq
\end{eqnarray*}
\begin{eqnarray*}
&\leq& \inf\{ t \mid \delta_n \varphi\to t \in T\}=|\delta_n \varphi|.
\end{eqnarray*}
Conversely,
\begin{eqnarray*}
|\delta_n \varphi|&=&\sup \{ t \mid  t \to \delta_n \varphi \in T\}
\leq \sup \{ t \mid  nt \to \varphi \in T\} \\
&=& \sup \{ \frac s{n} \mid  s\to \varphi \in T \} \\
&=& \frac{\sup \{ s \mid
   s \to \varphi \in T\}}{n}=\delta_n |\varphi|_T.
\end{eqnarray*}
\end{proof}

\begin{theorem} [Pavelka-style Completeness]
For R{\L}-theory $T$
$$
|\varphi|_T=||\varphi||_T.
$$
\end{theorem}
\begin{proof}
Soundness (i.e., $|\varphi|_T\leq ||\varphi||_T$) easily
follows from definition:
\begin{eqnarray*}
|\varphi|_T&=&\sup \{ r \mid r\leq e(\varphi)
\hbox{ with }e \hbox{ model of }T\}\\
&\leq & \inf \{ e(\varphi)\mid e \hbox{ model of }T\}=||\varphi||_T.
\end{eqnarray*}

Suppose without loss of generality
that $T$ is a consistent R{\L}-theory. Then there
exists a prime extension $T'\supseteq T$.
By Theorem \ref{evaluation},
the function $e: \varphi \in Form \to |\varphi|_T \in [0,1]$
is a model of $T'$,
and $|\varphi|_T=e(\varphi)\geq ||\varphi||_T$.
\end{proof}

\section{Complexity Issues}
In \cite{mun_NP} the SAT problem for {\L}ukasiewicz logic
is proved to be NP-complete.
In this section we shall prove that
the tautology
problem for Rational {\L}ukasiewicz logic is in co-NP
and since tautology problem of {\L}ukasiewicz
formulas can be reduced to tautology problem
of Rational {\L}ukasiewicz formulas as a subset, then
the latter is co-NP-complete.
Such result will be a byproduct of the fact that
if $\Gamma$ is a
finite set of {\L}ukasiewicz formulas and
$\varphi$ is a {\L}ukasiewicz formula, then
the problem to establish if $\Gamma \vdash \varphi$
is in co-NP  (see, for example, \cite{aguth}, \cite{agucia99}).
By \cite{wojc}, in this case
$\Gamma \vdash \varphi$ if and only if
for every assignment $v$ satisfying every formula
of $\Gamma$, $v(\varphi)=1$.

\medskip
Let us consider an
alphabet containing $\delta$ and a
symbol $|$ in such a way that
$\delta\underbrace{||\ldots |}_{n \hbox{ times}}$
stands for
$\delta_n$.

Let $\tau$ be a formula of Rational {\L}ukasiewicz logic,
with variables
among $\{X_1,$ $\ldots,$ $X_n\}$. Using
the same notation as in Subsection \ref{variety}, let
$\delta_{i_1}\tau_{i_1}, \ldots, \delta_{i_m}\tau_{i_m}$
denote all nodes of the parsing tree of subformulas of $\tau$
that begin with the symbol
$\delta$.

Let $\tau^*(X_1,\ldots, X_n,Z_1,\ldots,Z_m)$ be the formula
obtained by substituting every subformula $\delta_{i_j}\tau_{i_j}$
by the new variable $Z_j$.

Let $\Gamma$ be the set of {\L}ukasiewicz formulas defined by
$$
\Gamma = \bigcup_{j=1}^m\{i_j.Z_{i_j}\leftrightarrow \tau^*_{i_j},
\neg Z_{i_j} \odot (i_j-1).Z_{i_j} \},
$$
where $\tau^*_{i_j}$ has been obtained as $\tau^*$,
accordingly substituting occurrences of $\delta_{h_k}$ by
new variables $Z_{h_k}$. Then the formula $\tau$
is satisfiable in Rational {\L}ukasiewicz logic if and only if
$\Gamma \vdash \tau^* $ holds.
Since this last problem is in co-NP,
we have to give an estimation of lengths of $\Gamma$ and $\Delta$
in terms of the length of $\tau$.

\begin{definition}
\label{length}
The length of a formula of Rational {\L}u\-ka\-sie\-wi\-cz
logic is inductively defined
as follows:
\begin{itemize}
\item[(i)]
For every variable $X_i$, $\#X_i=1$
\item[(ii)]
$\#(\varphi\oplus \psi)=\#\varphi + \#\psi$
\item[(iii)]
$\#(\neg \varphi)=\#\varphi$
\item[(iv)]
$\#(\delta_n \varphi)=n+\#\varphi$
\end{itemize}
\end{definition}

Since this definition is an extension of the definition
of length for {\L}ukasiewicz formulas, we shall use the notation
$\#\varphi$ also when $\varphi$ is a {\L}ukasiewicz formula.
We set, without
loss of generality,
$\#(\varphi \leftrightarrow \psi)=
2(\#\varphi +\#\psi)$. If $\Lambda$ is a finite
set of formulas then
$$
\#\Lambda=\sum_{\lambda \in \Lambda}\#\lambda.
$$

If
$\delta_{i_1}\tau_{i_1},$ $\ldots,$ $\delta_{i_m}\tau_{i_m}$
are all subformulas of $\tau$ involving a connective $\delta$,
we have
\begin{equation}
\label{lengthtau*}
\#\tau^* \leq \#\tau -\sum_{j=1}^m i_m + m\leq \#\tau,
\end{equation}
because $\tau^*$ is obtained from $\tau$
by removing all occurrences of $\delta_{i_j}$.
The first inequality in (\ref{lengthtau*})
holds because not every new variable $Z_{i_j}$ appears
in $\tau^*$.
The second inequality holds since
$\sum_{j=1}^m i_m \geq m$.

In $\Gamma$ there is a pair of formulas
$i_j.Z_{i_j}\leftrightarrow \tau^*_{i_j},
\neg Z_{i_j} \odot (i_j-1).Z_{i_j}$
 for every $\delta_{i_j} \tau_{i_j}$ occurring
in $\tau$. Since
$$
\#(\delta_{i_j} \tau_{i_j})=i_j+\#\tau_{i_j}
$$
and
\begin{eqnarray*}
\#(i_j.Z_{i_j}\leftrightarrow \tau^*_{i_j})&=&2(i_j+\#(\tau^*_{i_j}))
\leq 2(i_j+\#(\tau_{i_j}))\\
\#(\neg Z_{i_j} \odot (i_j-1).Z_{i_j})&=&1+i_j-1,
\end{eqnarray*}
then
\begin{eqnarray*}
\label{gamma}
\# \Gamma&=&\sum_{j=1}^m \left( \# (i_j.Z_{i_j}
\leftrightarrow \tau^*_{i_j})
+\#(\neg Z_{i_j} \odot (i_j-1).Z_{i_j})\right)
\leq\\
&\leq &
\sum_{j=1}^m \left( 2(i_j+\#(\tau_{i_j}))+i_j\right)\leq
\sum_{j=1}^m 3\#(\delta_{i_j}\tau_{i_j})\\
&\leq & 3\#\tau.
\end{eqnarray*}
Putting together Equations (\ref{gamma}) and (\ref{lengthtau*})
we get the desired conclusion.

\bigskip
We shall now show that the complexity of the
tautology problem for Rational {\L}ukasiewicz logic
does not change if the index $n$ of $\delta_n$ is written
in binary notation.

Then let, in Definition \ref{length},
$\#\delta_n \varphi=\log_2 n+\#\varphi$.
If $\delta_n \tau_n$ occur in $\tau$
let $m_1>\ldots>m_h\geq 0$ be integer numbers
(depending on $n$) such that
$$
n=2^{m_1}+\ldots+2^{m_h}.
$$
We introduce $m_1$ new variables $Y_1,\ldots,Y_{m_1}$
and new formulas
\begin{eqnarray*}
\sigma(1,n)&=&Y_1\oplus Y_1\leftrightarrow  Y_2 \\
\ldots&&\\
\sigma(m_1-1,n)&=&Y_{m_1-1}\oplus
Y_{m_1-1}\leftrightarrow Y_{m_1}\\
\sigma'(1,n)&=&Y_1\odot Y_1\leftrightarrow0\\
\ldots&&\\
\sigma'(m_1-1,n)&=&Y_{m_1-1}\odot Y_{m_1-1}\leftrightarrow 0\\
\tau^*&=&\tau_n\leftrightarrow \left(2Y_{m_1}
\oplus 2Y_{m_2} \oplus \ldots \oplus
\left\{
\begin{array}{ll}
2Y_{m_h} &\hbox{ if }m_{h}>0\\
Y_{m_h} &\hbox{ if }m_{h}=0
\end{array}
\right. \right)
\end{eqnarray*}

We have $2(m_1-1)\leq 2 \log_2 n$ formulas $\sigma(i,n)$
and $\sigma'(i,n)$ of constant length and further
$\#\tau^*=2(\#\tau +2(m_1+\ldots+m_h))\leq 2(\#\tau +2\log_2 n)$.
Since $\#(\delta_n \tau)=\#\tau +\log_2 n$, then
$\#\tau^*\leq 4 \#(\delta_n \tau)$.

If
$\delta_{i_1}\tau_{i_1}, \ldots, \delta_{i_m}\tau_{i_m}$
are all subformulas of $\tau$ that begin with a symbol
$\delta$, then for any $\delta_{i_j}\tau_{i_j}$ we
suitably introduce formulas $\sigma,\sigma',\tau^*$
and thus reduce the problem of tautology
to the problem of deciding if a {\L}ukasiewicz formula
is consequence of a finite set of formulas.
The latter is co-NP
in the length of $\tau$.

%

\begin{thebibliography}{99}

\bibitem{aguth}
S.~Aguzzoli.
\newblock {\em Geometrical and Proof Theoretical Issues in {\L}ukasiewicz
  propositional logics}.
\newblock PhD thesis, 1998.
\newblock University of Siena, Italy.

\bibitem{agucia99}
S.~Aguzzoli and A.~Ciabattoni.
\newblock Finiteness in infinite-valued {{\L}}ukasiewicz logic.
\newblock {\em Journal of Logic, Language and Information}, 9:5--29, 2000.

\bibitem{agumun01}
S.~Aguzzoli and D.~Mundici.
\newblock Weierstrass approximations by {{\L}}ukasiewicz formulas with one
  quantified variable.
\newblock In {\em Proceedings of the 30th IEEE International Symposium on
  Multiple-Valued Logic (ISMVL)}. IEEE Computer Society Press,
pp. 361--366, 2001.

\bibitem{amapor}
P.~Amato and M.~Porto.
\newblock An algorithm for the automatic generation of a logical formula
  representing a control law.
\newblock {\em Neural Network World}, 10:777--786, 2000.

\bibitem{baazveith}
M.~Baaz and H.~Veith.
\newblock Interpolation in fuzzy logic.
\newblock {\em Archive for Mathematical logic}, 38:461--489, 1998.

\bibitem{chakie}
C.~C. Chang and H.~J. Keisler.
\newblock {\em Model theory}.
\newblock North-Holland Publishing Co., Amsterdam, 1973.

\bibitem{Chang}
C.C. Chang.
\newblock Algebraic analysis of many valued logics.
\newblock {\em Trans. Amer. Math. Soc.}, 88:467--490, 1958.

\bibitem{Chang2}
C.C. Chang.
\newblock A new proof of the completeness of the {{\L}}ukasiewicz axioms.
\newblock {\em Trans. Amer. Math. Soc.}, 93:74--90, 1959.

\bibitem{cigdotmun}
R.~Cignoli, I.M.L. {D'Ottaviano}, and D.~Mundici.
\newblock {\em Algebraic Foundations of many-valued reasoning}, volume~7 of
  {\em Trends in Logic}.
\newblock Kluwer, Dordrecht, 2000.

\bibitem{dinola}
A.~{Di Nola}.
\newblock Representation and reticulation by quotients of {MV}-algebras.
\newblock {\em Ricerche di Matematica}, 40:291--297, 1991.

\bibitem{dinlet96}
A.~{Di Nola} and A.~Lettieri.
\newblock Coproduct {MV}-algebras, nonstandard reals and {R}iesz spaces.
\newblock {\em Journal of Algebra}, 185:605--620, 1996.

\bibitem{Hajek}
P.~H\'{a}jek.
\newblock {\em Metamathematics of Fuzzy Logic}.
\newblock Trends in Logic. Kluwer, Dordrecht, 1998.

\bibitem{hoehle}
U.~H{\"o}hle.
\newblock Commutative residuated monoids.
\newblock In U.~H{\"o}hle and P.~Klement, editors, {\em Non-classical logic and
  their applications to fuzzy subsets}. Kluwer, 1995.

\bibitem{McNaughton}
R.~Mc{N}aughton.
\newblock A theorem about infinite-valued sentential logic.
\newblock {\em Journal of Symbolic Logic}, 16:1--13, 1951.

\bibitem{mun_interp}
D.~Mundici.
\newblock Interpretation of {AF} ${C}^*$-algebras in {{\L}}ukasiewicz
  sentential calculus.
\newblock {\em Journal of Functional Analysis}, 65:15--63, 1986.

\bibitem{mun_NP}
D.~Mundici.
\newblock Satisfiability in many-valued sentential logic is {NP}-complete.
\newblock {\em Theoretical Computer Science}, 52:145--153, 1987.

\bibitem{mun_constr}
D.~Mundici.
\newblock A constructive proof of {McNaughton}'s {Theorem} in infinite-valued
  logics.
\newblock {\em Journal of Symbolic Logic}, 59:596--602, 1994.

\bibitem{munpan}
D.~Mundici and G.~Panti.
\newblock Decidable and undecidable prime theories in infinite-valued logic.
\newblock {\em Annals of Pure and Applied Logic}, 108:269--278, 2001.

\bibitem{pavelka}
J.~Pavelka.
\newblock On fuzzy logic {I,II,III}.
\newblock {\em Zeitschr. f. Math. Logik und Grundl. der Math.},
  25:45--52,119--134,447--464, 1979.

\bibitem{wojc}
R.~W{\'o}jcicki.
\newblock On matrix representations of consequence operations of
{\L}ukasiewicz sentential calculi
\newblock {\em Zeitschr. f. Math. Logik und Grundl. der Math.},
19:239--247, 1973.
\newblock Reprinted in R. W{\'o}jcicki, G. Malinoski (Editors)
{\bf Selected Papers on {\L}ukasiewicz Sentential Calculi},
p. 101-111, 1977.
\end{thebibliography}
%

\end{document}